\def\autori{P.\ DALL'AGLIO}
\def\titolo{Stability results for solutions of obstacle problems with
measure data}

\font\sixrm=cmr6
\newcount\tagno \tagno=0                        
\newcount\thmno \thmno=0                        
\newcount\bibno \bibno=0                        
\newcount\chapno\chapno=0                       
\newcount\verno            
\newif\ifproofmode
\proofmodetrue
\newif\ifwanted
\wantedfalse
\newif\ifindexed
\indexedfalse

\def\ifundefined#1{\expandafter\ifx\csname+#1\endcsname\relax}
\def\Wanted#1{\ifundefined{#1} \wantedtrue 
\immediate\write0{Wanted #1  \the\chapno.\the\thmno}\fi}
\def\Increase#1{{\global\advance#1 by 1}}

\def\Assign#1#2{\immediate
\write1{\noexpand\expandafter\noexpand\def
 \noexpand\csname+#1\endcsname{#2}}\relax
 \global\expandafter\edef\csname+#1\endcsname{#2}}

\def\pAssign#1#2{\write1{\noexpand\expandafter\noexpand\def
 \noexpand\csname+#1\endcsname{#2}}}
\def\lPut#1{\ifproofmode\llap{\hbox{\sixrm #1\ \ \ }}\fi}
\def\rPut#1{\ifproofmode$^{\hbox{\sixrm #1}}$\fi}

\def\chp#1{\global\tagno=0\global\thmno=0\Increase\chapno
\Assign{#1}
{\the\chapno}{\lPut{#1}\the\chapno}}

\def\thm#1{\Increase\thmno 
\Assign{#1}
{\the\chapno.\the\thmno}\the\chapno.\the\thmno\rPut{#1}}

\def\frm#1{\Increase\tagno
\Assign{#1}{\the\chapno.\the\tagno}\lPut{#1}
{\the\chapno.\the\tagno}}

\def\bib#1{\Increase\bibno
\Assign{#1}{\the\bibno}\lPut{#1}{\the\bibno}}

\def\pgp#1{\pAssign{#1/}{\the\pageno}}

\def\ix#1#2#3{\pAssign{#2}{\the\pageno}
\immediate\write#1{\noexpand\idxitem{#3}
{\noexpand\csname+#2\endcsname}}}

\def\rf#1{\Wanted{#1}\csname+#1\endcsname\relax\rPut {#1}}

\def\rfp#1{\Wanted{#1}\csname+#1/\endcsname\relax\rPut{#1}}

\input \jobname.auxi
\Increase\verno
\immediate\openout1=\jobname.auxi

\immediate\write1{\noexpand\verno=\the\verno}

\ifindexed
\immediate\openout2=\jobname.idx
\immediate\openout3=\jobname.sym 
\fi


\font\twelverm=cmr12
\font\twelvei=cmmi12
\font\twelvesy=cmsy10
\font\twelvebf=cmbx12
\font\twelvett=cmtt12
\font\twelveit=cmti12
\font\twelvesl=cmsl12

\font\ninerm=cmr9
\font\ninei=cmmi9
\font\ninesy=cmsy9
\font\ninebf=cmbx9
\font\ninett=cmtt9
\font\nineit=cmti9
\font\ninesl=cmsl9

\font\eightrm=cmr8
\font\eighti=cmmi8
\font\eightsy=cmsy8
\font\eightbf=cmbx8
\font\eighttt=cmtt8
\font\eightit=cmti8
\font\eightsl=cmsl8

\font\sixrm=cmr6
\font\sixi=cmmi6
\font\sixsy=cmsy6
\font\sixbf=cmbx6

\catcode`@=11 
\newskip\ttglue

\def\twelvepoint{\def\rm{\fam0\twelverm}
\textfont0=\twelverm  \scriptfont0=\ninerm  
\scriptscriptfont0=\sevenrm
\textfont1=\twelvei  \scriptfont1=\ninei  \scriptscriptfont1=\seveni
\textfont2=\twelvesy  \scriptfont2=\ninesy  
\scriptscriptfont2=\sevensy
\textfont3=\tenex  \scriptfont3=\tenex  \scriptscriptfont3=\tenex
\textfont\itfam=\twelveit  \def\it{\fam\itfam\twelveit}%
\textfont\slfam=\twelvesl  \def\sl{\fam\slfam\twelvesl}%
\textfont\ttfam=\twelvett  \def\tt{\fam\ttfam\twelvett}%
\textfont\bffam=\twelvebf  \scriptfont\bffam=\ninebf
\scriptscriptfont\bffam=\sevenbf  \def\bf{\fam\bffam\twelvebf}%
\tt  \ttglue=.5em plus.25em minus.15em
\normalbaselineskip=15pt
\setbox\strutbox=\hbox{\vrule height10pt depth5pt width0pt}%
\let\sc=\tenrm  \let\big=\twelvebig  \normalbaselines\rm}

\def\tenpoint{\def\rm{\fam0\tenrm}
\textfont0=\tenrm  \scriptfont0=\sevenrm  \scriptscriptfont0=\fiverm
\textfont1=\teni  \scriptfont1=\seveni  \scriptscriptfont1=\fivei
\textfont2=\tensy  \scriptfont2=\sevensy  \scriptscriptfont2=\fivesy
\textfont3=\tenex  \scriptfont3=\tenex  \scriptscriptfont3=\tenex
\textfont\itfam=\tenit  \def\it{\fam\itfam\tenit}%
\textfont\slfam=\tensl  \def\sl{\fam\slfam\tensl}%
\textfont\ttfam=\tentt  \def\tt{\fam\ttfam\tentt}%
\textfont\bffam=\tenbf  \scriptfont\bffam=\sevenbf
\scriptscriptfont\bffam=\fivebf  \def\bf{\fam\bffam\tenbf}%
\tt  \ttglue=.5em plus.25em minus.15em
\normalbaselineskip=12pt
\setbox\strutbox=\hbox{\vrule height8.5pt depth3.5pt width0pt}%
\let\sc=\eightrm  \let\big=\tenbig  \normalbaselines\rm}

\def\ninepoint{\def\rm{\fam0\ninerm}
\textfont0=\ninerm  \scriptfont0=\sixrm  \scriptscriptfont0=\fiverm
\textfont1=\ninei  \scriptfont1=\sixi  \scriptscriptfont1=\fivei
\textfont2=\ninesy  \scriptfont2=\sixsy  \scriptscriptfont2=\fivesy
\textfont3=\tenex  \scriptfont3=\tenex  \scriptscriptfont3=\tenex
\textfont\itfam=\nineit  \def\it{\fam\itfam\nineit}%
\textfont\slfam=\ninesl  \def\sl{\fam\slfam\ninesl}%
\textfont\ttfam=\ninett  \def\tt{\fam\ttfam\ninett}%
\textfont\bffam=\ninebf  \scriptfont\bffam=\sixbf
\scriptscriptfont\bffam=\fivebf  \def\bf{\fam\bffam\ninebf}%
\tt  \ttglue=.5em plus.25em minus.15em
\normalbaselineskip=11pt
\setbox\strutbox=\hbox{\vrule height8pt depth3pt width0pt}%
\let\sc=\sevenrm  \let\big=\ninebig  \normalbaselines\rm}

\def\eightpoint{\def\rm{\fam0\eightrm}
\textfont0=\eightrm  \scriptfont0=\sixrm  \scriptscriptfont0=\fiverm
\textfont1=\eighti  \scriptfont1=\sixi  \scriptscriptfont1=\fivei
\textfont2=\eightsy  \scriptfont2=\sixsy  \scriptscriptfont2=\fivesy
\textfont3=\tenex  \scriptfont3=\tenex  \scriptscriptfont3=\tenex
\textfont\itfam=\eightit  \def\it{\fam\itfam\eightit}%
\textfont\slfam=\eightsl  \def\sl{\fam\slfam\eightsl}%
\textfont\ttfam=\eighttt  \def\tt{\fam\ttfam\eighttt}%
\textfont\bffam=\eightbf  \scriptfont\bffam=\sixbf
\scriptscriptfont\bffam=\fivebf  \def\bf{\fam\bffam\eightbf}%
\tt  \ttglue=.5em plus.25em minus.15em
\normalbaselineskip=9pt
\setbox\strutbox=\hbox{\vrule height7pt depth2pt width0pt}%
\let\sc=\sixrm  \let\big=\eightbig  \normalbaselines\rm}

\def\twelvebig#1{{\hbox{$\textfont0=\twelverm\textfont2=\twelvesy
	\left#1\vbox to10pt{}\right.\n@space$}}}
\def\tenbig#1{{\hbox{$\left#1\vbox to8.5pt{}\right.\n@space$}}}
\def\ninebig#1{{\hbox{$\textfont0=\tenrm\textfont2=\tensy
	\left#1\vbox to7.25pt{}\right.\n@space$}}}
\def\eightbig#1{{\hbox{$\textfont0=\ninerm\textfont2=\ninesy
	\left#1\vbox to6.5pt{}\right.\n@space$}}}
 
\def\displayliness#1{\null\,\vcenter{\openup1\jot \m@th
  \ialign{\strut\hfil$\displaystyle{##}$\hfil
    \crcr#1\crcr}}\,}
	       
\def\displaylinesno#1{\displ@y \tabskip=\centering
   \halign to\displaywidth{ \hfil$\@lign \displaystyle{##}$ \hfil
	\tabskip=\centering
     &\llap{$\@lign##$}\tabskip=0pt \crcr#1\crcr}}
		     
\def\ldisplaylinesno#1{\displ@y \tabskip=\centering
   \halign to\displaywidth{ \hfil$\@lign \displaystyle{##}$\hfil
	\tabskip=\centering
     &\kern-\displaywidth
     \rlap{$\@lign##$}\tabskip=\displaywidth \crcr#1\crcr}}

\catcode`@=12 

\def\parag#1#2{\goodbreak\bigskip\bigskip\noindent
		   {\bf #1.\ \ #2}
		   \nobreak\bigskip} 
\def\intro#1{\goodbreak\bigskip\bigskip\goodbreak\noindent
		   {\bf #1}\nobreak\bigskip\nobreak}
\long\def\th#1#2{\goodbreak\medskip\noindent
		{\bf Theorem #1.\ \ \it #2}\smallskip}
\long\def\lemma#1#2{\goodbreak\bigskip\noindent
		{\bf Lemma #1.\ \ \it #2}}
\long\def\prop#1#2{\goodbreak\bigskip\noindent
		  {\bf Proposition #1.\ \ \it #2}}

\long\def\defin#1#2{\goodbreak\bigskip\noindent
		  {\bf Definition #1.\ \ \rm #2}}
\long\def\rem#1#2{\goodbreak\medskip\noindent
		 {\bf Remark #1.\ \ \rm #2}}
\long\def\ex#1#2{\goodbreak\medskip\noindent
		 {\bf Example #1.\ \ \rm #2}}

\def\negbigskip{\vskip-\bigskipamount}

\def\proof{\vskip.4cm\noindent{\it Proof.\ \ }}
\def\proofof#1{
		\vskip.4cm\noindent{\it Proof of #1.\ \ }}

\def\sqr#1#2{\vbox{
   \hrule height .#2pt 
   \hbox{\vrule width .#2pt height #1pt \kern #1pt 
      \vrule width .#2pt}
   \hrule height .#2pt }}
\def\square{\sqr74}

\def\endproof{{\unskip\nobreak\hfill \penalty50
\hskip1em\hbox{}\nobreak\hfill $\square$ \goodbreak
\parfillskip=0pt  \finalhyphendemerits=0}}

\mathchardef\emptyset="001F
\mathchardef\hyphen="002D


\def\wto{\rightharpoonup}


\def\rightheadline{\eightpoint\hfil\titolo
\hfil\tenrm\folio} 
\def\leftheadline{\tenrm\folio\hfil\eightpoint
\autori \hfil}
\def\zeroheadline{\hfill} 
%
\headline={\ifnum\pageno=0 \zeroheadline
\else\ifodd\pageno\rightheadline
\else\leftheadline\fi\fi}

\nopagenumbers
\magnification=1200
\baselineskip=14pt
\hfuzz=2pt
\parindent=2em
\mathsurround=1pt
\tolerance=1000


\pageno=0
\hsize 14truecm
\vsize 25truecm
\hoffset=0.8truecm
\voffset=-1.55truecm

\null
\vskip 2.8truecm
{\twelvepoint
\baselineskip=1.7\baselineskip
\centerline{\bf STABILITY RESULTS FOR SOLUTIONS OF}
\centerline{\bf OBSTACLE PROBLEMS WITH MEASURE DATA}
}
\vskip2truecm
\centerline{Paolo DALL'AGLIO}
\vfil

{\eightpoint
\baselineskip=1.2\baselineskip
\centerline{\bf Abstract}
\bigskip
\noindent 
In this paper we study the continuous dependence with respect to
obstacles for obstacle problems with measure data. This is deeply
investigated introducing a suitable type of convergence, which gives
stability under very general hypotheses. Moreover stability with
respect to $\rm H^1$ and uniform convergent obstacles is proved.
\par
}

\vfil
\vskip 1truecm 
\centerline {Ref. S.I.S.S.A. 
142/99/M (December 1999)}
\vskip 1truecm 
\eject


\topskip=15pt 
\vsize 22.5truecm
\hsize 17truecm
\hoffset=0truecm
\voffset=0.5truecm

\overfullrule=0pt
\proofmodefalse 

\def\R{{\rm I\! R}}
\def\Rn{{\rm I\! R}^N}
\def\Luno{{\rm L}^1(\Omega)}
\def\Linf{{\rm L}^\infty(\Omega)}
\def\Hunozero{{\rm H}^1_0(\Omega)}   
\def\Huno{{\rm H}^1(\Omega)}
\def\Hduale{{\rm H}^{\hbox{\kern 1truept \rm -}\kern -1truept1}\kern
	-1truept(\Omega)}
\def\Wunoqzero{{\rm W}^{1,q}_0(\Omega)}
\def\Wunoq{{\rm W}^{1,q}(\Omega)} 
\def\Mb{{\cal M}_b(\Omega)}
\def\Mbp{{\cal M}^+_b(\Omega)}
\def\Mbo{{\cal M}^0_b(\Omega)}
\def\Mbop{{\cal M}^{0,+}_b(\Omega)}
\def\Inters{\Mb\cap\Hduale}
\def\intl{\int\limits}
\def\dualita#1#2{\langle#1,#2\rangle}

\def\Fgpsimu{{\cal F}^g_\psi(\mu)}  
\def\Ggpsimu{{\cal G}^g_\psi(\mu)}  
\def\wto{\rightharpoonup}

\def\cc{\subset\subset}
\def\um{u_{\mu}}
\def\un{u_{\nu}}
\def\ul{u_{\lambda}}
\def\uln{u_{\lambda_n}}
\def\ugo{u^g_0}
\def\uho{u^h_0}

\def\weak2{\hbox{\ weakly in }\Hunozero}
\def\strongq{\hbox{\ strongly in }\Wunoq}
\def\strong2{\hbox{\ strongly in }\Huno}
\def\weaks{\hbox{\ weakly-$\ast$ in }\Mb}
\def\mosco{\buildrel{\rm M}\over{-\!\!\!-\!\!\!\longrightarrow}}
\def\lev{\buildrel{\rm lev}\over{-\!\!\!-\!\!\!\longrightarrow}}
\def\ae{\hbox{\ a.e.\ in }\Omega}
\def\qe{\hbox{\ q.e.\ in }\Omega}

\def\vuoto{{\scriptstyle\bigcirc\!\!\!}\!\lower.15ex\hbox{/}}
\def\f{\varphi}
\def\e{\varepsilon}

\def\Kpsi{{\rm K}_{\psi}(\Omega)}
\def\Kpsin{{\rm K}_{\psi_n}(\Omega)}

\def\crr{\cr\noalign{\smallskip}}

\parag{\chp{intro}}{Introduction}

Given a regular bounded open set $\Omega$ of ${\R}^N$, $N\ge1$, and a
linear elliptic operator ${\cal A}$ of the form
$$
{\cal A}u =-\sum_{j,j=1}^N D_i ( a_{ij} D_j u ),
$$
with $a_{ij}\in L^\infty(\Omega)$, we study obstacle problems for the
operator ${\cal A}$ in
$\Omega$ with
homogeneous Dirichlet boundary conditions on $\partial\Omega$, when the
datum
$\mu$ is a bounded Radon measure on $\Omega$ and the obstacle
$\psi$ is an arbitrary function on $\Omega$.

According to~[\rf{DAL-LEO}], a function $u$ is a solution of this problem,
which will be denoted by $OP(\mu,\psi)$, if $u$ is the smallest
function with the following properties: $u\geq\psi$ in $\Omega$ and
$u$ is a solution in the sense of Stampacchia~[\rf{STA}] of a problem
of the form
$$
\cases{{\cal A}u=\mu+\lambda\quad&in $\Omega$\crr
	u=0&on $\partial\Omega$,\cr}\eqno(\frm{stampintro})
$$
for some bounded Radon measure $\lambda\geq 0$. The measure $\lambda$
which corresponds to the solution of the obstacle problem is called
the obstacle reaction.

Existence and uniqueness of the solution of $OP(\mu,\psi)$ have been
proved in~[\rf{DAL-LEO}], provided that there exists a measure
$\lambda$ such that the solution of~(\rf{stampintro}) is greater than
or equal to $\psi$. These results have been extended to the non-linear
case in~[\rf{LEO}], when $\mu$ vanishes on all sets with capacity
zero. For a different approach to obstacle problems for non-linear
operators with measure data see~[\rf{BOC-GAL}], [\rf{BOC-CIR-1}], 
[\rf{BOC-CIR-2}], [\rf{OPP-ROS-1}] and~[\rf{OPP-ROS-2}].

If the measure $\mu$ belongs to the dual $\Hduale$ of the Sobolev
space $\Hunozero$, and if there exists a function $w\in\Hunozero$
above the obstacle $\psi$, then the solution of the obstacle problem 
$OP(\mu,\psi)$ according to the previous definition coincides with the
solution $u$ of the variational inequality
$$
\cases{u\in\Hunozero,\ u\geq\psi,\crr
	\dualita{{\cal A}u}{v-u}\geq\dualita{\mu}{v-u},\crr
	\forall v\in\Hunozero\hbox{ s.t. }v\geq\psi\cr}
$$
where $\langle\cdot,\cdot\rangle$ denotes the duality pairing between
$\Hduale$ and $\Hunozero$. In this case the obstacle reaction
$\lambda$ belongs to $\Hduale$. We also know it is concentrated on the contact set
$\{u=\psi\}$ if $\psi$ is continuous, or, more in general, quasi upper
semicontinuous.

An important role in this problem is played by the space $\Mbo$ of all
bounded Radon measures on $\Omega$ which are absolutely continuous
with respect to the harmonic capacity. If the datum $\mu$ belongs to
$\Mbo$ (it is actually enough that its negative part $\mu^-$ is such), 
so does the obstacle reaction, provided that there exists  a
measure $\lambda\in\Mbo$ such that the solution of~(\rf{stampintro})
is greater than or equal to $\psi$ (see~[\rf{DAL-LEO}], theorem~7.5). In
this case the obstacle reaction is concentrated on the contact set
$\{u=\psi\}$, whenever the obstacle $\psi$ is quasi upper
semicontinuous (see~[\rf{LEO}], theorem~2.9). 

It can be seen that in general this does not occur when
$\mu^-\not\in\Mbo$. This case was studied
in~[\rf{DAL-DAL}]. Remembering that $\mu^-$ can be decomposed as
$\mu^-=\mu^-_a+\mu^-_s$, where $\mu^-_a \in{\cal M}_b^0(\Omega)$ and
$\mu^-_s$ is concentrated on a set of capacity zero, it was proved
that under some natural assumptions on the obstacle, the singular part
$\mu^-_s$ can be neglected and the obstacle problems $OP(\mu,\psi)$
and $OP(\mu^+-\mu^-_a,\psi)$ have the same solutions.

The topic of continuous dependence with respect to data was already
treated in~[\rf{DAL-LEO}].

As for stability with respect to the right hand side, it was proved
that, if $\mu_n,\,\mu\in\Mb$ are such that $\mu_n\to\mu$ strongly in
$\Mb$, then $u_n\to u$ strongly in $\Wunoq$, where $u_n$ and $u$ are
the solutions of $OP(\mu_n,\psi)$ and $OP(\mu,\psi)$ respectively.
Trying to use weak-$*$ convergence, it was seen that in general
$\mu_n\wto\mu$ weakly-$*$ does not imply that $u_n\to u$, even with
the obstacle $\psi\equiv 0$, but we know only that for any measure
$\mu\in\Mb$, there exists a special sequence $\mu_k\wto\mu$ weakly-$*$
in $\Mb$, with $\mu_k\in\Inters$, such that $u_k\to u$ strongly in
$\Wunoq$.

In this paper we consider stability with respect to obstacles. To
study this question we introduce a kind of convergence of functions,
the level set convergence, which yields the convergence of solutions
under very mild assumptions.

The convergence of $\psi_n$ to $\psi$ in the sense of level sets, defined
precisely in definition~\rf{livelli}, is verified in particular when
$$
{\rm cap}(\{\psi>t\}\cap B)=\lim_{n\to+\infty}{\rm cap}(\{\psi_n>t\}\cap
B)
$$
for  all $t\in\R$ and for all $B\cc \Omega$ (see also
remark~\rf{semplice}).

We will see that without further hypothesis it can only be proved that,
calling $u_n$ and $u$ the solutions of $OP(\mu,\psi_n)$ and of
$OP(\mu,\psi)$ respectively, if $\psi_n\lev\psi$ then, up to a
subsequence, $u_n$ converges to some function $u^*$ which is always
greater than or equal to $u$ (proposition~\rf{ustar}).

Then we will obtain, from the level set convergence of the obstacles,
that $u_n$ converges to $u$, under some conditions: in particular, by
means of the Mosco convergence of convex sets, we obtain that

\item{(i)} if $\mu^-\in\Hduale$ then $u_n\to u\strong2$;

\item{(ii)} if $\mu^-\in\Mbo$ then $u_n\to u\strongq$;

\item{(iii)} if $\psi$ is suitably controlled below, then $u_n\to u\strongq$.

In section~\rf{teoremini} we consider the case $\psi_n\leq\psi$ and
show that we have the convergence of solutions for any datum $\mu\in\Mb$.

We conclude this study considering two cases in which the assumptions
that the obstacles converge in a stronger way allows to obtain a
stronger convergence also for the solutions. When the difference
$\psi_n-\psi$ belongs to $\Hunozero$ and tends to zero strongly in
this space, then we obtain the same type of convergence for the
solutions, for any $\mu\in\Mb$.

In section~\rf{elleinfinito}, we extend the theory so
far developed to the case of nonzero boundary values. For any function
$g\in\Huno$, we can define the function $u$ to be the solution of
$OP(\mu,g,\psi)$ if and only if $u-\ugo$ is the solution of
$OP(\mu,\psi)$,
where $\ugo$ is the solution of
$$
\cases{	{\cal A}\ugo=0\quad&in $\Hduale$\crr
	\ugo-g\in\Hunozero.&\cr}
$$
All the results developed in the case of homogeneous boundary
conditions can be extended, thanks to the linearity of ${\cal A}$.

Using this extension we prove a new characterization: the solution of
$OP(\mu,g,\psi)$ is the minimum element among all the supersolutions of
${\cal A}-\mu$ which are above the obstacle and greater than or equal to
$g$ on the boundary $\partial\Omega$. 
From this we easily prove that if the obstacles 
converge uniformly then so do the solutions of the corresponding
obstacle problems for any $\mu\in\Mb$.

\parag{\chp{fattinoti}}{Notations and basic results.}
				 
Let $\Omega$ be an open bounded subset of $\Rn$, $N\geq1$, with
Lipschitz boundary.

Let ${\cal A}(u)=-{\rm div}(A(x)\nabla u)$ be a linear elliptic operator
with coefficients in $L^\infty(\Omega)$, that is $A(x)=(a_{ij}(x))$ is
an $N\times N$ matrix such that
$$
a_{ij}\in L^\infty(\Omega)\hbox{\ \ and \ }\sum
a_{ij}(x)\xi_i\xi_j\geq\gamma
|\xi|^2,\quad\forall\xi\in\R^N, \ae,
$$
with $\gamma>0$.

We want to consider the obstacle problem also in the case of thin
obstacles, so we will need the techniques of capacity theory. For this
theory we refer, for instance, to [\rf{HEI-KIL}].

We recall very briefly that, given a set $E\subseteq\Omega$ its
capacity with respect to $\Omega$ is given~by
$$
{\rm cap}(E)=\inf\{\|z\|^2_{\Huno}\,:\,z\in\Hunozero, z\geq 1\hbox{ a.e.\ in a
neighbourhood of }E\}.
$$

A property holds quasi everywhere (abbreviated as q.e.) when it holds up
to sets of capacity zero. 

A set $A$ is said to be quasi open (resp.\ quasi closed) if for any $\e>0$
there exists an open (resp.\ closed) set $V$ such that cap$(A\triangle
V)<\e$.

A function $v:\Omega\to{\overline{\R}}$ is quasi continuous (resp.
quasi upper semicontinuous) if, for every $\e>0$ there exists a set   
$E$ such that ${\rm cap}(E)<\e$ and $v|_{\Omega\setminus E}$ is continuous
(resp. upper semicontinuous) in
$\Omega\setminus E$.

We recall also that if $u$ and $v$ are quasi continuous functions and
$u\leq v$ a.e.\ then also $u\leq v$ q.e.\ in $\Omega$.

A function $u\in\Hunozero$ always has a quasi continuous
representative, that is there exists a quasi continuous function
$\tilde u$ which equals $u$ a.e.

Consider the function $\psi:\Omega\to{\overline{\R}}$, and let the convex
set be
$$
K_{\psi}(\Omega):=\{z\ \hbox{quasi continuous }  : z\geq \psi \hbox{
q.e.\ in }\Omega\}.
$$
  
Without loss of generality we may suppose that $\psi$ is quasi upper
semicontinuous thanks to the following proposition (it is a
consequence of proposition 1.5 in [\rf{DAL-4}]).

\prop{\thm{invilup}}{Let $\psi:\Omega\to{\overline{\R}}$.  Then 
there exists a quasi upper semicontinuous function $\hat
\psi:\Omega\to\overline{\R}$ such that:

\item{ 1.} $\hat \psi\geq\psi$ q.e.\ in $\Omega$;

\item{ 2.} if $\varphi:\Omega\to{\overline{\R}}$ is quasi upper
semicontinuous and $\varphi\geq\psi$ q.e.\ in $\Omega$ then
$\varphi\geq\hat\psi$ q.e.\ in~$\Omega$.}

Thus, in particular, $K_{\psi}(\Omega)=K_{\hat{\psi}}(\Omega)$.

In their natural setting, obstacle problems are part of
the theory of variational
inequalities (for which we refer to well known books such as
[\rf{KIN-STA}] and [\rf{TRO}]).

For any datum $F\in\Hduale$ the
variational inequality with obstacle $\psi$
$$
\cases{\dualita{{\cal A} u}{v-u}\geq\dualita{F}{v-u}
		\quad\forall v\in\Kpsi\cap\Hunozero\crr
	u\in\Kpsi\cap\Hunozero\cr}\eqno(\frm{vi})
$$
(which, for simplicity, will be indicated by
$
VI(F,\psi)
$),
has a unique solution, whenever the set $\Kpsi\cap \Hunozero$ is nonempty,
i.e.
$$
\exists z\in\Hunozero\ :\ z\geq\psi\quad\qe.\eqno(\frm{zeta})  
$$
In this case we will say that the obstacle is $VI$-admissible.

Let now $\Mb$ be the space of bounded Radon measures, with the norm
given by $\|\mu\|_{\Mb}=|\mu|(\Omega)$. $\Mbo$ is the subspace of
measures of $\Mb$ vanishing on sets of zero capacity. $\Mbp$ and
$\Mbop$ are the corresponding cones of non negative measures.  Recall
that $\Hduale\not\subseteq\Mb$ but $\Hduale\cap\Mb\subseteq\Mbo$.
  
Any measure $\mu\in\Mb$ can be decomposed as $\mu=\mu_a+\mu_s$
(see~[\rf{FUK}]), where $\mu_a\in\Mbo$ and $\mu_s$ is concentrated on
a set of capacity zero.

If $x\in\Omega$, we denote by $\delta_x$ the Dirac's delta centered at
$x$.

When the datum is a measure, equations and inequalities can not be studied
in the variational framework.

In [\rf{STA}] G. Stampacchia gave the following definition

\defin{\thm{stamp}}{}A function
$\um\in\Luno$ is a solution
in the sense of Stampacchia (also called
solution by duality) of the
equation
$$
\cases{{\cal A}\um=\mu\ \ & in $\Omega$\crr
	\um=0& on $\partial\Omega$\cr}\eqno(\frm{stampacchia})
$$
if
$$
\intl_\Omega \um g\,dx\,=\,\intl_\Omega u^*_g\,d\mu,\quad\forall
g\in\Linf,
$$
where $u^*_g$ is the solution of
$$
\cases{{\cal A}^*u^*_g=g\quad \hbox{ in }\Hduale\crr
	u^*_g\in\Hunozero\cr}
$$
and ${\cal A}^*$ is the adjoint of ${\cal A}$.

Throughout the paper $q$ will be any exponent satisfying
$1<q<{N\over{N-1}}$; Stampacchia proved that a solution
$\um$ exists, is unique, and belongs to $\Wunoqzero$; moreover if the
datum $\mu$ is more regular, namely belongs to $\Inters$, then the
solution coincides with the variational one. It is possible to prove
that, when the data converge weakly-$*$ in $\Mb$, the solutions
converge strongly in $\Wunoqzero$.
  
If $u$ is such a solution then $T_k(u)\in\Hunozero$, for any
$k\in\R^+$, where $T_k(s):=(-k)\vee(s\wedge k)$ is the usual
truncation function. Moreover
$$
\intl_\Omega|DT_k(u)|^2dx\leq k\,|\mu|(\Omega).\eqno(\frm{tronche})
$$
These facts imply that $u$ has a quasi continuous representative which is
finite$\qe$. In the rest of the paper we shall always identify $u$ with
its quasi continuous representative.

We will use the following notation: $\um$ denotes
the solution of the equation~(\rf{stampacchia}).

The following definition of solution for obstacle problems with measure
data was given in [\rf{DAL-LEO}].

\defin{\thm{disvar}}{}We say that the function
$u$ is a solution of the
obstacle problem with datum $\mu$ and obstacle $\psi$ (shortly
$OP(\mu,\psi)$) if
\item{1.} $u\in\Kpsi$ and there exists a positive bounded measure
$\lambda\in\Mbp$ such
that
$$
u=\um+\ul;
$$
\item{2.} for any $\nu\in\Mbp$, such that $v=\um+\un$ belongs to $\Kpsi$,
we have
$$
u\leq v\ae.
$$

For the problem to make sense let us assume that the obstacle $\psi$
satisfies a minimal
hypothesis, instead of (\rf{zeta}), namely
$$
\exists \rho\in\Mb\ :\ u_\rho\geq\psi\hbox{ q.e.\ in
}\Omega.\eqno(\frm{ipomin})
$$
In this case we will say that $\psi$ is $OP$-admissible.

\th{\thm{esist2}}{Given $\mu\in\Mb$ and $\psi$ $OP$-admissible
(namely $\psi\leq u_\rho$ q.e.\ in $\Omega$), the
obstacle problem $OP(\mu,\psi)$ has a unique solution $u=\um+\ul$.
The measure $\lambda\in\Mbp$ satisfies
$$
\|\lambda\|_{\Mb}\leq\left\|(\mu-\rho)^-\right\|_{\Mb}\eqno(\frm{mumenorho}).
$$
Moreover, if $\rho\in\Inters$, there exists a sequence $\mu_k={\cal
A}T_k(\um-u_\rho)+\rho
\in\Inters$ such that the solutions
$u_k$ of
$OP(\mu_k,\psi)$ converge strongly in $\Wunoq$. This sequence does not
depend on the obstacle but only on the measure $\rho$.}

Also here the positive measure $\lambda$ associated with the solution is
called the obstacle reaction.

We mention here a very simple and very useful result whose proof is
immediate, but it is worth stating it on its own.

\lemma{\thm{magico}}{Let $\mu,\,\nu\in\Mb$ and $\psi$ is $OP$-admissible. 
Then $u$ is the solution of $OP(\mu+\nu,\psi)$ if and only if $u-\un$ is
the solution of $OP(\mu,\psi-\un)$}.
  
The particular case in which the datum $\mu$ belongs to $\Mbo$ has
been investigated in detail. Then the obstacle reaction itself belongs
to $\Mbo$, provided the obstacle satisfies the following condition
$$
\exists \sigma\in\Mbo\ :\ u_\sigma\geq\psi\hbox{ q.e.\ in
}\Omega;\eqno(\frm{ipominzero})
$$
This will be shortened by saying that $\psi$ is $OP^0$-admissible.

Notice that if the datum $\mu$ is in $\Mbo$, but the obstacle is
only $OP$-admissible, then the reaction $\lambda$ in general does
not belong to $\Mbo$. This is shown in the following simple example.

\ex{\thm{lambda}}{Let $\mu=0$ and $\psi=u_{\delta_{x_0}}$, where
$\delta_{x_0}$ is the Dirac's delta centered at $x_0\in\Omega$. Then the
solution of $OP(0,\psi)$ is $u_{\delta_{x_0}}$ itself and hence
$\lambda=\delta_{x_0}\not\in\Mbo$.}

The interaction between obstacles and solutions is studied deeply in
[\rf{DAL-DAL}], where the following result was proved.

\th{\thm{lostesso}}{Let $\mu\in\Mb$ and let $\psi:\Omega\to\overline\R$ be
such that
$$
-u_{\tau}-u_{\sigma}-\f\leq\psi\leq u_{\sigma}\quad\qe\eqno(\frm{condizione})
$$
where $\f\in\Huno$, $\sigma\in\Mbo$ and $\tau\in\Mb$ such that
$\tau\perp\mu^-_s$.
Then the
solutions 
$$
u=\um+\ul \hbox{ of } OP(\mu,\psi)\ \ \ \hbox{ and }\ \ \ 
u_{\mu^+}-u_{\mu^-_a}+u_{\lambda_1}\hbox{ of }OP(\mu^+-\mu^-_a,\psi)
$$
are the same.  Moreover $\lambda=\lambda_1+
\mu^-_s$ with $\lambda_1\in\Mbop$.}

Here condition (\rf{condizione}) is given in its full generality, for
instance it is
satisfied by obstacles in $\Huno$ that are $OP$-admissible.

\rem{\thm{delte}}{If the obstacle $\psi$ satisfies~(\rf{condizione})
with $\tau=0$ then the conclusion holds for every $\mu\in\Mb$.

The presence of $\tau$, which deppends on $\mu$, in (\rf{condizione})
allows anyway to treat situations like the following one.  If ${\cal
A}=-\Delta$, $\Omega=B_1(0)$, the obstacle is $-u_{\delta_0}$ and the
datum is $-\delta_{x_0}$ for any $x_0\neq 0$, then the solution of the
obstacle problem is zero, because the theorem applies, and because the
solution must be less than or equal to zero.  }

In this paper we will be concerned with the continuous dependence of
the solutions with respect to various types of convergence of the
obstacles. To our knowledge the only result that was proved on this
problem for an arbitrary measure $\mu$ is the following (proved
in~[\rf{DAL-LEO}]) which deals with a very special case.

\prop{\thm{psienne}}{Let $\psi_n:\Omega\to\overline{\R}$ be obstacles
such
that
$$
\psi_n\leq\psi \quad\hbox{ and }\quad\psi_n\to\psi\quad\qe,
$$
$\psi$ $OP$-admissible, and let $u_n$ and $u$ be the solutions of
$OP(\mu,\psi_n)$ and $OP(\mu,\psi)$, respectively. Then
$$
u_n\to u\quad\hbox{ strongly in }\Wunoq.
$$}
{\vskip-1.5em}

This result will be generalized in section~\rf{teoremini}.

In order to study the problem in the most general way we introduce in section~\rf{mosco}
a notion of convergence that was used in~[\rf{DAL-4}]. It will be called
``level set convergence'' and in most cases it is equivalent to the
convergence of convex sets introduced by U.  Mosco in~[\rf{MOS}], see
also~[\rf{ATT}]

Mosco proved that this type of convergence is the right one for the
stability of variational inequalities with respect to obstacles. This is
the main theorem of his theory.

\th{\thm{vimosco}}{Let $\psi_n$ and $\psi$ be $VI$-admissible. Then
$$
\Kpsin\cap\Hunozero\mosco\Kpsi\cap\Hunozero,
$$
if and only if, for any $f\in\Hduale$,
$$
u_n\to u\strong2
$$
where $u_n$ and $u$ are the solutions of $VI(f,\psi_n)$  and
$VI(f,\psi)$, respectively.}

\parag{\chp{mosco}}{The level set convergence, and the related stability
properties.}

In this section we will define a kind of convergence of functions
which will prove to be a good one for the obstacles in obstacle
problems with measure data: it is rather general and allows to obtain
the convergence of the solutions under very mild assumptions.

\defin{\thm{livelli}}{Let $\psi_n$ and $\psi$ be quasi upper semicontinuous
function from $\Omega$ to $\overline\R$. We say that $\psi_n$ tends to
$\psi$ in the sense of level sets and write
$$
\psi_n\lev\psi
$$
if 
$$
\quad\qquad {\rm cap}(\{\psi>t\}\cap B)\leq\liminf_{n\to+\infty}
{\rm cap}(\{\psi_n>s\}\cap B')\eqno(\frm{liminf})
$$
$$
{\hskip-1.5cm}\limsup_{n\to+\infty}{\rm cap}(\{\psi_n>t\}\cap B)\leq
{\rm cap}(\{\psi>s\}\cap B')\eqno(\frm{limsup})
$$
for all $s,t\in\R$, $s<t$, and for all $B\subset\subset B'\subset
\subset \Omega$.}

\rem{\thm{semplice}}{From the definition it is clear that the level set
convergence is implied~by 
$$
{\rm cap}(\{\psi>t\}\cap B)=\lim_{n\to+\infty}{\rm
cap}(\{\psi_n>t\}\cap B).
\eqno(\frm{limitecap})
$$
for all $t\in\R$ and for all $B\cc\Omega$}
\rem{\thm{cap}}{From the definition it follows that, if
$\psi_n$ converge to $\psi$ locally in capacity,~i.e.
$$
{\rm cap}(\{|\psi_n-\psi|>t\}\cap A)\to0,\quad\forall t\in\R^+,\;
\forall A\cc\Omega,
$$
then $\psi_n\lev\psi$.}

To prove that the limit in the sense of level set is unique we give the
following lemma from capacity theory which can be found in~[\rf{FUG}].
\lemma{\thm{capaci}}{Let $E$ and $F$ be quasi closed subsets of $\Omega$ 
such that
$$
{\rm cap}(E\cap A)\leq {\rm cap}(F\cap A),\quad\forall A\subseteq\Omega\hbox{ open },
\eqno(\frm{cap})
$$
then ${\rm cap}(E\setminus F)=0$ (we say also that $E$ is quasi contained
in $F$).}

\prop{\thm{benposto}}{Let $\psi_n$, $\psi$ and $\f$ be quasi upper
semicontinuous functions. If
$$
\psi_n\lev\psi\qquad\hbox{ and }\qquad\psi_n\lev\f,
$$
then $\psi=\f$.}
\proof 
Let us fix an open set $A\cc\Omega$ and two real numbers  $s<t$. Take now
two subsets $A'$ and $A''$ such that $A''\cc A'\cc A$ and real numbers
$t'$ and $t''$ such that $s<t'<t''<t$. Then
$$
\eqalign{
{\rm cap}(\{\psi>t''\}\cap A'')&
	\leq \liminf_{n\to+\infty}{\rm cap}(\{\psi_n>t'\}\cap A')\crr
&\leq\limsup_{n\to+\infty}{\rm cap}(\{\psi_n>t'\}\cap A')
		\leq {\rm cap}(\{\f\geq s\}\cap A).\cr}
$$
Hence, since $\{\psi\geq t\}\subseteq\{\psi>t''\}$, we have
$$
{\rm cap}(\{\psi\geq t\}\cap A'')\leq {\rm cap}(\{\f\geq s\}\cap A),
$$
from which, invading $A$ by means of $A''\cc A$,
$$
{\rm cap}(\{\psi\geq t\}\cap A)\leq {\rm cap}(\{\f\geq s\}\cap A).
$$
Using the fact that $\psi$ and $\f$ are quasi upper semicontinuous
and thanks to lemma~\rf{capaci} we deduce that $\{\psi\geq t\}$ is quasi
contained $\{\f\geq s\}$. Now, given $t$, consider two sequences
$t_k\searrow t$ and $s_k\searrow t$, with $t_k> s_k$, so that 
$$
\{\psi\geq t_k\}\nearrow \{\psi> t\}\hbox{\ \ and\ \ }
\{\f\geq s_k\}\nearrow \{\f> t\},
$$
and we get that 
$$
\{\psi> t\}\hbox{ is quasi contained in }\{\f>t\},
$$
for all $t\in\R$.

Exchanging the roles of $\psi$ and $\f$ we get the reverse inclusion so
that
$\{\psi> s\}$ and $\{\f>s\}$ coincide up to sets of capacity zero.

Now we recover the
values of $\psi$ and $\f$ at quasi every point $x\in\Omega$ thanks to the
well known formula
$$
\f(x)=\sup_{s\in Q}s\,\chi_{\{\f> s\}}(x).
$$
Since the level sets are the same, the two functions coincide quasi
everywhere.
\endproof

The main result on level sets convergence is the following theorem,
which shows the connection with the Mosco convergence introduced
in~[\rf{MOS}] (for the proof see theorem~5.9 in~[\rf{DAL-4}]).

\th{\thm{equiva}}{Let 
$\psi_n$ and $\psi$ be functions $\Omega\to\overline\R$.
If  
$$
\Kpsin\cap\Hunozero\mosco\Kpsi\cap\Hunozero.
$$ 
then 
$$
\psi_n\lev\psi.
$$
If moreover the obstacles are equicontrolled
from above, namely
$$
\psi_n,\ \psi\,\leq u_\rho\,,
\quad\hbox{ with }\rho\in\Inters,
$$
then also the reverse implication holds.}

Notice that, though very similar to Mosco convergence, the level set
convergence concerns also the case of obstacles that are not
$VI$-admissible.

Another simple observation, which requires no proof, but which is
useful to state separately is the following.
\lemma{\thm{crescente}}{Let $\psi_n\lev\psi$ and let $\Phi:\R\to\R$ be a
continuous non decreasing function. Then
$$
\Phi(\psi_n)\lev\Phi(\psi).
$$}
{\vskip-1.5em}

In the next lemmas we will denote the solution of $OP(\mu,\psi_n)$ and
of $OP(\mu,\psi)$ by $u_n$ and $u$, respectively.

Let us show that in general the Mosco convergence (and so also the level
set convergence) of the obstacles does not imply the convergence of the
solutions for an arbitrary measure.

\ex{\thm{nomosco}} Let $\Omega=B_1(0)\subseteq\R^N$, 
with $N>2$, ${\cal A}=-\Delta$ and $\mu= 
-\delta_0$, the Dirac delta in the origin.

Let the obstacles $\psi_n=-n$, so that clearly 
$$
\Kpsin\cap\Hunozero\mosco\Hunozero
$$
and, by theorem~\rf{equiva}, also $\psi_n\lev -\infty$

It is immediate to see that the solutions $u_n=u_{-\delta_0}
+u_{\lambda_n}$ of $OP(-\delta_0, -n)$ are less than or equal to zero
since the latter satisfies condition 1 of definition~\rf{disvar}.  So
$u_n=T_n(u_n)$ and hence is in $\Hunozero$. But then
$-\delta_0+\lambda_n\in\Hduale\cap\Mb\subset\Mbo$, and it must be a
positive measure and hence $u_n=0$ for each $n$.  On the other hand
$u=u_{-\delta_0}$ and cannot be the limit of the~$u_n$.

What can be proved without further assumptions is the following result.

\prop{\thm{ustar}}{Let $\psi_n,\psi\leq u_\rho\qe$ with 
$\rho\in\Mb$. Assume that 
$$
\psi_n\lev\psi
$$
Then there exists a subsequence $u_{n'}$ and a quasi continuous function
$u^*\in\Wunoqzero$, such that
$$
u_{n'}\to u^\ast\strongq,
$$
and 
$$
u^\ast\geq u\quad\qe.
$$}
{\vskip-1cm}
\proof By theorem~\rf{esist2}
$$
\|\lambda_n\|_{\Mb}\leq\left\|(\mu-\rho)^-\right\|_{\Mb}\eqno(\frm{lambdan})
$$
so that there exists a subsequence
$\{\lambda_{n'}\}$ and a measure
$\lambda^*\in\Mbp$ such that $\lambda_n\wto\lambda^\ast,
\weaks$ and hence $u_{n'}=\um+u_{\lambda_{n'}}\to
u^\ast=\um+u_{\lambda^\ast}$ strongly in $\Wunoq$.

If we show that $u^\ast\geq\psi\qe$ we will have that $u^\ast$
satisfies condition~1 of definition~\rf{disvar} and get the thesis,
by definition~\rf{disvar}.

Given $k>0$, observe that, thanks to (\rf{tronche}) and to (\rf{lambdan})
$$
 \intl_\Omega|DT_k(u_{n'})|^2dx\leq kc
$$ 
and hence $T_k(u_{n'})\wto T_k(u^\ast)\weak2$.

From lemma~\rf{crescente} it follows that
$$
T_k(\psi_n)\lev T_k(\psi),\quad\forall k\in\R^+;
$$
Since $T_k(u_n)\geq T_k(\psi_n)\qe$ for each $n$ and $k$, and using
theorem~\rf{equiva} and 
the definition of Mosco convergence,  we get
$T_k(u^\ast)\geq T_k(\psi)\qe$.

Now we can pass to the limit as $k\to+\infty$ and obtain
$u^\ast\geq \psi\qe$.
\endproof

We prove now the central lemma of this section.

\lemma{\thm{psimenow}}{Let $\psi_n$ and $\psi$ be quasi upper
semicontinuous functions
controlled above by $u_\rho$ with $\rho\in\Inters$, and let $w$ be a quasi
continuous function. If 
$$
\psi_n\lev\psi
$$
then
$$
\psi_n-w\lev\psi-w.
$$}
To prove this lemma we need the following result.
 
\lemma{\thm{c2approx}}{Given a quasi continuous function $w$, for 
each $A\subset\subset\Omega$ and for each $\varepsilon>0$, there exists
$u\in C^{\infty}_0(\Omega)$ such that
$$
{\rm cap}\left({\left\{|u-w|>\varepsilon\right\}\cap A}\right)
<\varepsilon.\eqno(\frm{umenow})
$$}
{\vskip -1 truecm}
\proof By definition of quasi continuity there exists a relatively
closed subset $C$ such that 
${\rm cap}(\Omega\setminus C)<\varepsilon$ and $w|_{C}$  is
continuous. By Tietze's theorem there exists a continuous function
$g$ which extends $w_{|_{C\cap\overline A}}$ it to $\Rn$.

Obviously, for any $A\cc\Omega$, we have 
$\left\{|w-g|>0\right\}\cap A\subseteq\Omega\setminus C$ so
that
$${\rm cap}\left({\left\{|w-g|>0\right\}}\cap A\right)
<\varepsilon
$$. 
{\vskip-1.5em}

On its turn $g$ can be approximated in $A$ with a  
function $u\in C^\infty_0(\Omega)$ so
that $\sup_A|u-g|<\varepsilon$, and again from
$\left\{|w-u|>\varepsilon\right\}\cap A\subseteq
\left\{|w-g|>0\right\}\cap A$ we get 
${\rm cap}\left({\left\{|w-u|>\varepsilon\right\}}\cap A\right)<
\varepsilon$.
\endproof

\proofof{\rf{psimenow}}
It is immediate to observe that the thesis is true when, instead of 
$w$, we have a function $u$ in $\Hunozero$.
This is because in our hypotheses level set convergence is equivalent to
Mosco convergence (see theorem~\rf{equiva}) and translating $\Kpsin$ and
$\Kpsi$ by
$u\in\Hunozero$ we get
$$
{\rm K}_{(\psi_n-u)}\cap\Hunozero\mosco{\rm K}_{(\psi-u)}\cap\Hunozero.
$$
or equivalently
$$
\psi_n-u\lev\psi-u
$$

We want to show the inequalities (\rf{liminf}) and (\rf{limsup}) for
$\psi_n-w$ and $\psi-w$.

Let us fix $B\cc \Omega$, $\e>0$ and
a function $u\in C^{\infty}_0(\Omega)$ such that (\rf{umenow}) holds with
respect to
$B$.

Observe now that for any $t\in\R$ we have
$$
\left\{\psi_n-w>t\right\}\cap B\subseteq\left({\left\{
\psi_n-u>t-\e\right\}\cap B}\right)\cup
\left({\left\{u-w>\e\right\}\cap B}\right),
$$
hence, by subadditivity,
$$
{\rm cap}\left({\left\{\psi_n-w>t\right\}\cap B}\right)
\leq {\rm cap}\left({\left\{
\psi_n-u>t-\e\right\}\cap B}\right)+{\rm cap}
\left({\left\{|u-w|>\e\right\}\cap B}\right).
$$
Passing to the limsup and using (\rf{umenow}), we obtain
$$
	\limsup_{n\to+\infty} 
	{\rm cap}\left({\left\{\psi_n-w>t\right\}\cap B}\right)
	\leq \limsup_{n\to+\infty} {\rm cap}\left({\left\{
	\psi_n-u>t-\e\right\}\cap B}\right)+\e.
$$
We know that for $\psi-u$,~(\rf{liminf}) holds true,
so we can use it with $t-\e$ and $t-2\e$
instead of
$t$ and $s$, so that, for $B\cc B'\cc\Omega$, we get
$$
\limsup_{n\to\infty} 
{\rm cap}\left({\left\{\psi_n-w>t\right\}\cap B}\right)\leq
{\rm cap}\left({\left\{
\psi-u>t-2\e\right\}\cap B'}\right)+\e.\eqno(\frm{passaggio})
$$
With an argument similar to before, from
$$
\left\{\psi-u>t-2\e\right\}\cap B'\subseteq\left({\left\{
\psi-u>t-3\e\right\}\cap B'}\right)\cup
\left({\left\{w-u>\e\right\}\cap B'}\right),
$$
we obtain,
$$
{\rm cap}\left({\left\{\psi-u>t-2\e\right\}\cap B'}\right)
\leq {\rm cap}\left({\left\{
\psi-u>t-3\e\right\}\cap B'}\right)+\e,
$$
and substituting in (\rf{passaggio})
we get
$$
\limsup_{n\to\infty} 
{\rm cap}\left({\left\{\psi_n-w>t\right\}\cap B}\right)\leq
{\rm cap}\left({\left\{
\psi-w>t-3\e\right\}\cap B'}\right)+2\e.
$$
For any choice of $s$ and $t$, $\e$ can be taken sufficiently small so
that $s<t-3\e$. Then we can let $\e\to0$ and conclude
$$
\limsup_{n\to\infty} 
{\rm cap}\left({\left\{\psi_n-w>t\right\}\cap B}\right)\leq
{\rm cap}\left({\left\{
\psi-w>s\right\}\cap B'''}\right).
$$
Here nothing depends on $u$ so this holds for all $s,t\in\R$, $s<t$ and
for all $B\cc B'\cc\Omega$.

Inequality (\rf{liminf}) is proved in a similar way and this concludes the 
proof.
\endproof

\th{\thm{psimenoumu}}{Let $\psi_n,\psi\leq u_\rho$ with $\rho\in
\Mb\cap\Hduale$. Let $\mu\in\Mb$ with $\mu^-\in\Hduale$ 
and let $u_n$ and $u$ be the solutions of $OP(\mu,\psi_n)$
and of $OP(\mu,\psi)$ respectively.
If
$$
\psi_n\lev\psi
$$
then $u_n-u$ belongs to $\Hunozero$ and tends to zero 
strongly in $\Hunozero$.}
\proof
Thanks to the previous lemma we have that
$$
\psi_n-\um\lev\psi-\um,
$$
and then, since $\psi_n-\um,\,\psi-\um\leq 
u_{\rho}+u_{\mu^-}$ q.e.\ in
$\Omega$, by theorem~\rf{equiva}, 
$$
{\rm K}_{(\psi_n-\um)}\cap\Hunozero\mosco{\rm
K}_{(\psi-\um)}\cap\Hunozero.
$$
Hence all solutions of variational inequalities converge. In particular
if $v_n$ and $v$ are the solutions of $VI(0,\psi_n-\um)$ and 
$VI(0,\psi-\um)$, respectively, then $v_n\to v$ strongly $\Hunozero$. 

By lemma~\rf{magico} we have $u_n=v_n+\um$ and $u=v+\um$. This implies
that $u_n-u=v_n-v$ and the conclusion follows.
\endproof

The minimal hypothesis on the obstacles $\psi_n$ and $\psi$ in order to
have the solutions of $OP(\mu,\psi_n)$ and $OP(\mu,\psi)$ is that they are
$OP$-admissible. Nevertheless if in this theorem we drop the
request that they are controlled by a function which is also in
$\Hunozero$, the conclusion fails. Ideed there is the following example
which derives from example~\rf{nomosco}.

\ex{\thm{controllo}}{Let us consider the operator ${\cal A}=-\Delta$, the
domain $\Omega=B_1(0)\subseteq\R^N$, with $N>2$, 
and the datum $\mu=0$. 

Consider now as obstacles $\psi_n=u_{\delta_0}-n$, where ${\delta_0}$
is the Dirac delta centred at zero. They are
clearly
$OP$-admissible, and also bounded by the same function $u_\rho$, but in
this case $\rho=\delta_0\not\in\Hduale$.

Now for each $n$ the solution $u_n$ of $OP(0,\psi_n)$ is $u_{\delta_0}$
itself. Indeed, according to lemma~\rf{magico}, $u_n-u_{\delta_0}$ is the
solution of $OP(-\delta_0,-n)$, that, as seen in 
example~\rf{nomosco}, is zero.

But then we have that $\psi_n\lev -\infty$ and
$u_n\to u_{\delta_0}$, while the solution of $OP(0,-\infty)$ is $u=0$.
}

When the negative part of the measure $\mu$ is only in $\Mbo$, we can not
use the same trick because the the sets $K_{\psi_n-\um}$ might be empty,
but anyway we do not fall into the pathology of the example \rf{nomosco}, and
in fact we
can prove the following theorem which gives the convergence of the
solutions as well, though in a weaker sense.

\th{\thm{menoemmezero}}{Let $\psi_n,\psi\leq u_\rho$ with $\rho\in
\Mb\cap\Hduale$, $\mu\in\Mb$ such that $\mu^-\in\Mbo$, 
and let $u_n$ and $u$ be the solutions of $OP(\mu,\psi_n)$
and of $OP(\mu,\psi)$ respectively.
If
$$
\psi_n\lev\psi,
$$
then $u_n\to u\strongq$.}
\proof
From [\rf{DAL-3}] we know that $\mu^-$ can be written as $g\nu$ where
$\nu\in\Mbp\cap\Hduale$ and $g\in
L^1(\Omega,\nu)$, $g\geq 0$.
Hence the measure $\mu^-_k:=(g\wedge k)\nu$ is in $\Hduale$, so that
$\mu_k:=\mu^+-\mu^-_k$ satisfies the hypothesis of the previous theorem.

Call $u^k_n$ and $u^k$ the solutions of $OP(\mu^k,\psi_n)$ and 
$OP(\mu^k,\psi)$, respectively. By theorem~\rf{psimenoumu},
$$
u^k_n\to u^k,\ \strong2,\quad\forall k>0.
$$

Now, observing that $\mu^--\mu^-_k$ is a positive 
measure, we easily obtain by comparison 
that $u^k_n\geq u_n$ and $u^k\geq u,\qe$.

On the other hand
$$
u_n+u_{(\mu^--\mu^-_k)}=\um+u_{(\lambda_n+\mu^--\mu^-_k)},
\eqno(\frm{quella})
$$
where $\lambda_n\geq 0$ is the obstacle reaction of $OP(\mu,\psi_n)$.
Since also $u_n+u_{(\mu^--\mu^-_k)}\geq u_n\geq\psi_n$, by definition
\rf{disvar} and eq.~(\rf{quella}) we have $u_n+u_{(\mu^--\mu^-_k)}
\geq u^k_n\qe$. In the same way we prove that
$u+u_{(\mu^--\mu^-_k)}\geq u^k\qe$.

Since $\mu^--\mu^-_k\to 0$ strongly in $\Mb$, we have, by
proposition~4.2 in~[\rf{DAL-LEO}] $u_{(\mu^--\mu^-_k)}\to 0\strongq$,
so, from
$$
u+u_{(\mu^--\mu^-_k)}\geq u^k\geq u\quad\qe
$$
letting $k\to\infty$ we get that $u^k\to u\ae$.

Recalling proposition~\rf{ustar}, let us fix a subsequence $\{u_{n'}\}$
which converges to a function $u^*$ strongly in $\Wunoq$, so that from
$$
u_{n'}+u_{(\mu^--\mu^-_k)}\geq u^k_{n'} \geq u_{n'}\quad\qe
$$
letting first $n'\to\infty$ and then $k\to\infty$ we obtain $u^\ast
\geq u\geq u^\ast\qe$. Therefore $u^k_n\to u$, since the limit does not
depend on the subsequence.
\endproof

As seen in example \rf{controllo} the request that the obstacles be well
controlled can not be dropped, even if the datum is regular. On the other
hand example \rf{nomosco} showed that the control from above can be not
enough to have convergence for all data $\mu\in\Mb$.

In the following theorem we show how, provided we strengthen the
assumptions on the obstacles in the way given by theorem~\rf{lostesso},
we can give up any assumption on the data $\mu$.

Notice that in the examples is always the limit obstacle the one that
gives troubles.  Indeed we see here that it is enough to require the
control from below only for the limit.

\th{\thm{tutto}}{Let $\psi_n,\psi\leq u_\rho$ with $\rho\in
\Mb\cap\Hduale$, let $\mu\in\Mb$  and let $\psi$be an obstacle.
Let $u_n$ and $u$ be the solutions of $OP(\mu,\psi_n)$
and of $OP(\mu,\psi)$, respectively.
If $\psi$ and $\mu$ satisfy condition~(\rf{condizione}) and
\vskip-1em
$$
\psi_n\lev\psi
$$
\noindent then $u_n\to u\strongq$.}
\proof
From proposition~\rf{ustar} we know that, up to a subsequence,  $u_n\to
u^*$ strongly in $\Wunoq$, and
$u^*\geq u\qe$.

Now consider the $v_n$'s solutions of $OP(\mu^+-\mu^-_a,\psi_n)$. These, by
theorem
\rf{menoemmezero}, converge to $v$, the solution of
$OP(\mu^+-\mu^-_a,\psi)$, but, according to theorem~\rf{lostesso}, $v=u$.

On the other side $v_n=\um+\uln+u_{\mu^-_s}$, with $\lambda_n\in\Mbp$, and
$v_n\geq\psi_n\qe$ and so, by definition~\rf{disvar}, we have
$$
v_n\geq u_n\quad\qe.
$$
Letting $n$ go to $+\infty$ we obtain $u\geq u^*\qe$. Therefore $u_n\to
u\strongq$. The limit doesn't depend on the subsequence, so the whole
sequence $u_n$ converges to~$u$.
\endproof

Let us show now a further example, which clarifies more
deeply in which cases there is not convergence of the solutions.

In particular we see that theorems~\rf{lostesso} and~\rf{tutto} do not
hold for some obstacles $\psi$ which are too singular only at one point.

\ex{\thm{menodelta}}{}Let us choose ${\cal A}=-\Delta$,
$\Omega=B_1(0)\subseteq\R^N$, with $N>2$, and $\mu=-\delta_0$. Let us
consider the obstacles $\psi=-u_{\delta_0}$ and
$$
\psi_n(x)=\cases{- \displaystyle{1\over 2}
		u_{\delta_0}(x)\quad& if \ \ \ 
			$|x|<a_n$\crr
		-n& if \ $a_n<|x|<b_n$\crr
		-u_{\delta_0}& if \ $b_n<|x|$\cr}
\eqno(\frm{biz})$$
where $a_n$ and $b_n$ are appropriate constants, which tend to zero as
$n\to+\infty$ (see picture).

It is easy to verify that  $\psi_n\lev\psi$ and that the solution of the
limit problem $OP(-\delta_0,\psi)$ is clearly
$-u_{\delta_0}$  itself.

Let us prove that the solution of $OP(-\delta_0,\psi_n)$ 
is $- {1\over 2}u_{\delta_0}$.

This function satisfies condition~1 of definition~\rf{disvar} because it is of the form $\um+u_{{1\over
2}\delta_0}$ and it is above the obstacle for each $n$.

Fix $n$ and suppose $\nu_n\in\Mbp$ such that $\um+u_{\nu_n}$ is the
solution.
Then it is smaller than or equal to $\um+u_{{1\over 2}\delta_0}$, or also
$$
\un\leq u_{{1\over 2}\delta_0}\quad\qe.
$$
\eject
\hbox{\hsize 17cm
\vtop{\hsize 8truecm
In the small circle $B_{a_n}(0)$ they must be equal. In
$B_1(0)\setminus B_{a_n}(0)$, $\un$ is superharmonic and $u_{{1\over
2}\delta_0}$ is harmonic and they have the same boundary data. So
$\un\geq u_{{1\over 2}\delta_0}$ and they must coincide.

This proves that the solution $u_n$ of $OP(-\delta_0,\psi_n)$ is
$- {1\over 2}u_{\delta_0}$ independently of $n$,
and that $u_n$ does not converge
to the solution $u=u_{\delta_0}$ of $OP(-\delta_0,\psi)$. 

As of remark~\rf{delte} we point out that in the example it is crucial
that the
deltas involved are centered in the same point. If for instance, with the
same obstacles, we had as datum $\mu=-\delta_{x_0}$ for any $x_0\neq 0$, we
would obtain, thanks to theorem~\rf{lostesso} that the solutions of
$OP(-\delta_{x_0},\psi_n)$ and of $OP(-\delta_{x_0},\psi)$ are all
identically zero.

The last consideration of this section concerns the fact that passing
from theorem~\rf{psimenoumu} to theorem~\rf{menoemmezero} we loose
something on the convergence of the solutions. To see that this loss
is not due to the technique of the proof we can consider the following
example.}
\qquad
\vtop{\hsize 7truecm \vskip 11cm
\includegraphics{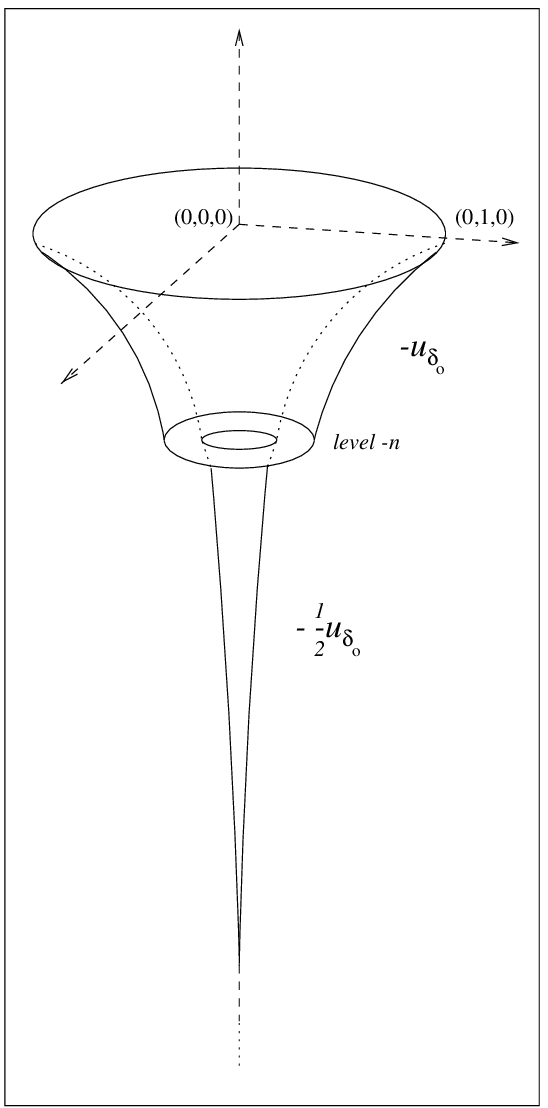}}}

\ex{\thm{orsina}}{}
Let $\Omega=B_{1\over 2}(0)\subseteq\R^N$, with $N>2$, 
and let $f\in\Luno$ be the function defined
by
$$
f(x)=\cases{
	\displaystyle{1\over{|x|^N(-\log|x|)^\vartheta}}\quad& if $x\in
		\Omega,\, x\neq 0$\crr
	0& if $x=0$\cr}
$$
with $\vartheta>1$. L. Orsina in [\rf{ORS}] noticed that the solution
$u_f$ of the equation 
$$
\cases{-\Delta u_f=f\quad &in $\Omega$\crr
	u_f=0&on $\partial\Omega$\cr }
$$
\vskip-.5em
\noindent belongs to $\Wunoq$ for any $q< {N\over N-1}$ but does not
belong  to ${\rm W}^{1,{N\over N-1}}(\Omega)$.

With this choice of ${\cal A}$ and $\Omega$, take as datum $\mu=-f$,
which clearly belongs to $\Mbo$, and as limit obstacle $\psi=-u_f$, which
satisfies
condition (\rf{condizione}) with $\sigma=f$, so that the solution of
$OP(-f,-u_f)$
is $-u_f$ itself.

If we set $\psi_n=-(u_f\wedge n)$ then the solution $u_n$ of
$OP(-f,-(u_f\wedge n))$ is between $0$ (because $f$ is positive and
$u_{-f}+u_f$ is a supersolution) and $-n$. Hence $u_n=T_n(u_n)$ and
this implies that $u_n\in\Hunozero$.

Now it is easy to see (use for instance remark~\rf{semplice}) that
$\psi_n\lev\psi$ so that, by proposition~\rf{tutto}, $u_n\to
u\strongq$. Nevertheless it is not possible to have this convergence
in the norm of ${\rm W}^{1,p}(\Omega)$ with $p\geq {N\over
N-1}$, because the fact that $u_n\in\Huno$ would imply also that
$u\in{\rm W}^{1,p}(\Omega)$, which is false since $u=-u_f$.

\vskip-1em

\parag{\chp{teoremini}}{The case $\psi_n\leq\psi$.}
\vskip-.5em

In the previous section we have seen that the level sets convergence
of the obstacles in general it is not enough to give the convergence
of the solutions for any $\mu$.

In this section we want to generalize
proposition~\rf{psienne}. Pointwise convergence is replaced by level
sets convergence, the condition $\psi_n\leq\psi\quad\qe$ is strong enough
to give the result for any measure $\mu\in\Mb$ with no control from below on
the limit obstacle.

\vskip-.5em

\prop{\thm{dasotto}}{Let $\psi_n$ and $\psi$ be $OP$-admissible and such 
that
\vskip-.5em
$$
\psi_n\lev\psi.
$$
\vskip-.5em
\noindent If in addition $\psi_n\leq\psi\qe$ then 
\vskip-.5em
$$
u_n\to u\quad\qe.
$$}
\vskip-1.5em
\proof
First apply proposition~\rf{ustar} from which we know that, up to a
subsequence, $u_n\to u^*\strongq$ and that $u^*\geq u$.  
The reverse inequality is guaranteed by the fact that $u\geq u_n$, for
all $n$ thanks to condition~1 of definition~\rf{disvar}.
\endproof

Let us remark that this is a generalized version of proposition
\rf{psienne}. Indeed under the assumption that $\psi_n\leq\psi$, we have
that quasi everywhere convergence implies level sets convergence. 

If $\psi_n$ is monotone increasing, this is an easy consequence of the
continuity of capacities on increasing sequences of sets.

The general case, $\psi_n\leq\psi$ but $\psi_n$ not necessarily
increasing, is proved by considering the sequence
$\f_n:=\inf_{k\geq n}\psi_k$, which is increasing, converge pointwise
to $\psi$ and satisfies
$$
{\rm cap}(\{\f_n>t\}\cap B)\leq
{\rm cap}(\{\psi_n>t\}\cap B)\leq
{\rm cap}(\{\psi>t\}\cap B)
$$
for all $t\in\R$, and $B\cc\Omega$. Passing to the limit, thanks to the
previous step, we conclude the proof.

\parag{\chp{convhuno}} {Obstacles converging in the energy space.}

It is well known that in the case of variational inequalities
the convergence of obstacles in the norm of $\Huno$ implies the 
convergence of the corresponding solutions. In particular we have the
following result

\th{\thm{vihuno}}{Let $\psi_1,\psi_2:\Omega\to\overline{\R}$ 
be 
$VI$-admissible. Suppose that $\psi_1-\psi_2\in\Hunozero$. Let $f\in\Hduale$  and let 
$u_1$ and $u_2$ be the solutions of $VI(f,\psi_1)$ and $VI(f,\psi_2)$, 
respectively. Then
$$
\|u_n-u\|_{\Hunozero}\leq {C\over\gamma}\,\left\|\psi_n-\psi\right\|_{\Hunozero},\eqno(\frm{acca})
$$
where $\gamma$ is the ellipticity constant and $C$ is such that
$|\dualita{{\cal A}u}{v}|\leq C\|u\|_{\Hunozero}\|v\|_{\Hunozero}$, so they  
depend only on the operator ${\cal A}$.} 

Also for the solutions of obstacle problems, we want to investigate
the dependence on $\Huno$ convergence. In this frame, as we have seen
with example~\rf{orsina}, this can not follow directly from Mosco
convergence, as it was in the variational case.  The next theorem
concerns the case in which the obstacles ``have the same boundary
value''.

\th{\thm{ophunozero}}{Let $\psi_n,\psi:\Omega\to
\overline{\R}$ be such that $\psi_n,\psi\leq u_\rho$ 
q.e.\ in $\Omega$, with $\rho\in\Mb$. 
Assume  $\psi_n-\psi\in\Hunozero$, and let $\psi_n-\psi
\to 0\strong2$. Let $\mu\in\Mb$ and let 
$u_n$ and $u$ be the solutions of $OP(\mu,\psi_n)$ and 
$OP(\mu,\psi)$, respectively. 
Then
$$
u_n-u\in\Hunozero\hbox{ \ and \ }u_n-u\to 0\strong2,
$$
}
\proof
As a first step assume that $\rho$ belongs also to $\Hduale$,
and consider the special sequence $\mu_k$ of 
measures in $\Mb\cap\Hduale$ 
$\mu_k\wto\mu,\weaks$, such that the solutions of 
the corresponding obstacle
problems converge (see theorem~\rf{esist2}). In particular $u^k_n\to
u_n,\strongq$ for all
$n$ and $u^k\to u,\strongq$.

Thanks to (\rf{acca}), for all $k$ we also have 
$$
\|u^k_n-u^k\|_{\Hunozero} \leq c\|\psi_n-\psi\|_{\Hunozero},
$$
so that the sequence $\{u_n^k-u^k\}_k$ is bounded in $\Hunozero$, 
for each $n$ fixed. Thus, up to a subsequence, there is a limit 
function $z$. But we already know that the sequence converges,
strongly in $\Wunoq$, to $u_n-u$, so this must be also the weak limit
in $\Huno$.

By lower semicontinuity of the norm we have
$$
\|u_n-u\|_{\Hunozero} \leq\liminf_{k\to\infty}
\|u^k_n-u^k\|_{\Hunozero} \leq c\|\psi_n-\psi\|_{\Hunozero}.
$$
This says that $u_n-u$ belongs to $\Hunozero$ (while $u_n$ and $u$, in 
general, do not) and gives the thesis in the first case.

Let now $\rho$ be only in $\Mb$.
Set
$\psi^h:=\psi\wedge h$ and $\psi^h_n:=\psi_n-\psi+\psi^h$. 
These obstacles are equi $OP$-admissible, because $\psi_n^h\leq\psi_n$ and
$\psi^h\leq\psi$ q.e.\ in $\Omega$.
They are also equi  $VI$-admissible since, if $\psi\leq u_\rho$, then 
$\psi^h\leq T_h(u_\rho)\in\Hunozero$ and $\psi^h_n\leq\psi_n-\psi+
T_h(u_\rho)\in\Hunozero$, and we can find a 
function $\f\in\Hunozero$ such
that $\f\geq\psi_n-\psi$ for all $n$. 

It is easy to see that in this case there 
exists $\rho_h\in\Inters$ such that $\psi^h_n,\,\psi^h\leq u_{\rho_h}$
q.e.\ in $\Omega$. Hence we are in the hypothesis of the
first step. 
Moreover $\|\psi^h_n-\psi^h\|_{\Hunozero}=\|\psi_n-\psi\|_{\Hunozero}$.

So by the first step, for each $h$
$$
\|u^h_n-u^h\|_{\Hunozero} \leq\gamma\|\psi_n-\psi\|_{\Hunozero}.
$$
On the other side we know that, since $\psi_n^h\nearrow\psi_n$ and
$\psi^h\nearrow\psi$, by proposition~\rf{psienne}  $u^h_n-u^h\to
u_n-u\strongq$, so that
we can conclude as in the first step.
\endproof

We want to remark that if, more generally, the obstacles are such that
$\psi_n-\psi\to 0$ in $\Huno$ then they also converge in the sense of
level sets, so we can deduce the convergence of the solutions in all the
situations given by theorems \rf{psimenoumu}, \rf{menoemmezero} and
\rf{tutto}, but here we obtain a stronger convergence with no further
assumptions on the obstacles and on the data.

We may now wonder what happens when the obstacles converge in the space
$\Wunoq$, with $1<q< {N\over N-1}$. In general this is not
enough to obtain the convergence of the
solutions. Indeed reconsider example~\rf{menodelta}. 
Let us prove that
$\psi_n\to\psi$ strongly in $\Wunoqzero$. We have
$$
\psi_n-\psi=\cases{ {\displaystyle{1\over 2}} u_{\delta_0} 
		& in $\quad|x|<a_n$\crr
		u_{\delta_0}-n\qquad& in $a_n<|x|<b_n$\crr
		0& otherwise\cr}
$$
so that 
$$
\|\psi_n-\psi\|^q_{\Wunoq}=
{1\over 2}\intl_{B_{a_n}(0)}|Du_{\delta_0}|^qdx+
{\hskip-.6cm}\intl_{B_{b_n}(0)\setminus
B_{a_n}(0)}{\hskip-.6cm}|Du_{\delta_0}|^qdx,
$$
which tends to zero, since $a_n$ and $b_n$ tend to zero and by
the absolute continuity of the integral.
But, as already seen in example~\rf{menodelta},
the solutions of the obstacle problems do not converge.

Anyway it is possible to prove the following result.

\prop{\thm{rhoenne}}{Let $\mu$ be in $\Mb$ and let $\psi_n$ and $\psi$ be 
$OP$-admissible and such that $\psi_n-\psi=u_{\rho_n}$ with
$\rho_n\in\Mb$, $\|\rho_n\|_{\Mb}\to 0$. Then
$$
u_n\to u\strongq,
$$
where $u_n$ and $u$ are the solutions of $OP(\mu,\psi_n)$ and
$OP(\mu,\psi)$, respectively.}
\proof
Since $\psi_n=\psi-u_{\rho_n}$, we have (using lemma \rf{magico}) that
$u_n-u_{\rho_n}$ is the solution of $OP(\mu+\rho_n,\psi)$. So from
theorem~4.2 in~[\rf{DAL-LEO}] we get that 
$$
u_n-u_{\rho_n}\to u\strongq.\eqno(\frm{urhon})
$$

Then
$$
\|u_n-u\|_{\Wunoq}\leq
\|u_n-u_{\rho_n}-u\|_{\Wunoq}+\|u_{\rho_n}\|_{\Wunoq};
$$
the first term goes to zero because of~(\rf{urhon}), 
the second one by hypothesis, and we get the thesis.
\endproof

\parag{\chp{elleinfinito}}{Problems with nonzero boundary data and
uniform convergence}

In this section we extend the theory of obstacle problems with measure
data developed in [\rf{DAL-LEO}] to problems with nonzero 
boundary data. This is standard for variational inequalities
and also in this case 
this generalization is very simple; we will
only point out what has to be settled.

Let $g\in\Huno$ we will denote by $u^g_0$ the solution of
$$
\cases{{\cal A}u^g_0=0\qquad& in $\Hduale$\crr
	u^g_0-g\in\Hunozero.\cr}
$$
We will look for solutions of obstacle problems which take the value $g$
on the boundary $\partial \Omega$. So we have to change accordingly the
notion of admissibility for the obstacles.

An obstacle
$\psi:\Omega\to{\overline{\R}}$ is said to be $OP_g$-admissible if
$$
\exists \rho\in\Mbp\ \hbox{s.t. }\psi\leq u_\rho+u^g_0\quad\qe.
$$

Given a measure $\mu\in\Mb$, a boundary datum $g\in\Huno$ and an
$OP_g$-admissible obstcale $\psi$, the solution of the obstacle
problem $OP(\mu,g,\psi)$, if it exists, is the minimum element of the
set
$$
{\cal F}^g_\psi(\mu):=\left\{v\in\Wunoq\,:\,\exists\nu\in\Mbp,\,
v=\um+u^g_0+\un;\,v\geq\psi\qe\right\}.
$$

It is immediate to prove the following

\th{\thm{esistg}}{Let $\mu\in\Mb$ and let $\psi$ be $OP_g$-admissible.
Then there exists a unique solution of $OP(\mu,g,\psi)$.}
\proof
Consider the obstacle $\psi-u^g_0$. It is $OP$-admissible. So there
exists a unique solution $v$ of $OP(\mu,\psi-u^g_0)$. Then $v+\ugo$ is our
solution: indeed it belongs to $\Fgpsimu$, and it is less than or equal to
any $z\in\Fgpsimu$.
\endproof

\rem{\thm{stima}}{From theorem~\rf{esistg} and (\rf{mumenorho}) it
follows also that, if $\um+\ugo+\ul$ is the solution of $OP(\mu,g,\psi)$,
then
$$
\|\lambda\|_{\Mb}\leq\left\|(\mu-\rho)^-\right\|_{\Mb},
$$
independently of $g$.}

\rem{\thm{convergenze}}{Since $\psi$ is $OP$-admissible if and only if
$\psi-\ugo$ is $OP_g$-admissible and
$u_n\to u\strongq$ if and only if $u_n-\ugo\to u-\ugo\strongq$, all
the theorems on continuous dependence on the data hold 
without modifications, in particular
propositions \rf{psienne} which will be useful in the 
following.}

\rem{\thm{coincide}}{When $\mu\in\Inters$ and $\psi\leq u_\rho+u^g_0\qe$
with $\rho\in\Inters$ then the solution of $OP(\mu,g,\psi)$ coincides
with the solution of $VI(\mu,g,\psi)$\footnote{$^\ast$}{\eightpoint There is no need to
define explicitly what is $VI(\mu,g,\psi)$.}.}

We come now to discuss the continuous dependence of the solutions on the
obstacles when these converge uniformly.

To do this we will use a characterization via supersolutions similar to
the one that holds in the variational case (see~[\rf{KIN-STA}]). 

%
%
%
%

To this aim, let us introduce the set ${\cal G}^g_\psi(\mu)$ of all
the functions $v\in\Wunoq$ with $v\geq\psi\qe$, such that
$v=\um+u^h_0+\un$, where $\nu\in\Mbp$ and $h\in\Huno$ such that $h\geq
g$ on $\partial\Omega$, i.e.\ $(h-g)^-\in\Hunozero$.

We see now that the solution of $OP(\mu,g,\psi)$ can be compared not only
with the functions of $\Fgpsimu$, but also with all those
that have boundary datum greater than or equal to $g$.

\prop{\thm{Ggpsimu}}{Let $\mu\in\Mb$ and $\psi$ be $OP_g$-admissible.
If $u$ is the solution of $OP(\mu,g,\psi)$ then it is the
minimum element of $\Ggpsimu$.}
\proof
{\it Step 1.} First consider $\mu\in\Inters$ and $\psi$ both
$VI$- and $OP$-admissible.

Let $v=\um+\uho+\un\in\Ggpsimu$. We approximate $\nu$ by means of the
sequence $\nu_k:={\cal A}T_k(\un)$. We have that $\nu_k\in\Inters$ and
that $\nu_k\wto\nu\weaks$. Moreover observe that $u_{\nu_k}=T_k(\un)$ tends
to $\un$ q.e.\ in $\Omega$ and, since $\un$ is nonnegative it is an
increasing sequence.

Hence if we define $v_k:=\um+\uho+u_{\nu_k}$, then $v_k\nearrow v\qe$, and
setting $\psi_k:=\psi\wedge v_k$ also $\psi_k\nearrow \psi\qe$

Let now $u_k$ be the
solutions of $VI(\mu,g,\psi_k)$. So, by theorem~II.6.4 in~[\rf{KIN-STA}],
$v_k\geq u_k$.
Using proposition~\rf{psienne} and remark~\rf{convergenze} 
we know that $u_k\to u\ae$. Then $v\geq
u\ae$ and then also$\qe$.

{\it Step 2.} Consider now $\mu\in\Mb$ and $\psi$ still 
both $VI$- and $OP$-admissible. Take again $v\in\Ggpsimu$.

Let $\mu_k={\cal A}T_k(\um-u_\rho)+\rho$ be  the sequence of measures 
given in theorem~\rf{esist2}, so that we know that if $u_k$ are the 
solutions of $VI(\mu_k,g,\psi)$ then $u_k\to u\strongq$. 

Taking now $v_k=u_{\mu_k}+\uho+\un$ it is easy to verify that $v_k\geq
\psi\qe$ for all $k>0$, and then, by definition~\rf{disvar}, $v_k
\geq u_k\qe$.
Also $v_k\to v\strongq$ so, passing to the limit, we obtain $v\geq u\ae$ 
and then also $\qe$.

{\it Step 3.} Finally consider the general case $\mu\in\Mb$ and $\psi$ 
$OP$-admissible.
The obstacles $\psi_k:=\psi\wedge k$ are also $VI$-admissible and such
that $\psi_k\nearrow\psi\qe$. So,  if 
$u_k$ is the solution of $OP(\mu,g,\psi_k)$, 
by proposition~\rf{psienne} and remark~\rf{convergenze}, 
we have that $u_k\to u\strongq$.

Taken any $v\in\Ggpsimu$, then $v\geq \psi_k$, for all $k$. Hence, by 
definition~\rf{disvar}, $v\geq u_k\qe$. Passing to the limit,
we get $v\geq u\ae$ and also$\qe$.
\endproof
From this we point out that the sets $\Fgpsimu$ and $\Ggpsimu$ have the
following lattice property.
 
\prop{\thm{reticolo}}{Let $\mu\in\Mb$, $g\in\Huno$ and $\psi$
$OP_g$-admissible. Then
\item{(i)}If $u,v\in\Fgpsimu$ then $u\wedge v\in\Fgpsimu$;
\item{(ii)}If $u,v\in\Ggpsimu$ then $u\wedge v\in\Ggpsimu$.}
\proof
Let us prove only the first statement, the proof of the second being
alike.

Set $w:=u\wedge v$ and let $z$ be the solution of $OP(\mu,g,w)$. Then
$u,v\in {\cal F}^g_w(\mu)$ and hence also $w\geq z$.

On the other hand $z\geq w$ and hence they are equal. So $u\wedge v$ is of
the form $\um+\ugo+\un$ and is above $\psi$, and hence belongs to
$\Fgpsimu$.
\endproof

We can prove now the following continuity result

\th{\thm{uniforme}}{Let $\mu\in\Mb$, $g\in\Huno$, and $\psi_n$ and
$\psi$ be
$OP$-admissible 
and let $u_n$ and $u$ be the solutions of $OP(\mu,\psi_n)$ and $OP(\mu,
\psi)$, respectively.
Assume that $\psi_n-\psi\in\Linf$ and
$\psi_n-\psi\to0$ in $\Linf$.
Then
$$
u_n-u\in\Linf\quad\hbox{ and }\quad u_n-u\to0\hbox{ in }\Linf.
$$}
\negbigskip
\proof
Set $c_n:=\|\psi_n-\psi\|_{\Linf}$. Obviously $c_n=u^{c_n}_0$, so that
$$
u+c_n=\um+u^{c_n+g}_0+\ul\ \hbox{ and }\ u+c_n\geq\psi_n\quad\qe
$$
hence $u+c_n\in{\cal G}^g_{\psi_n}(\mu)$ and hence $u+c_n\geq u_n$.

The same can be done the other way round to obtain that $u_n+c_n\geq u$.
In the end we get $|u_n- u|\leq c_n$, and, taking the sup over $x\in\Omega$, the
thesis.
\endproof

\rem{\thm{implica}}{Also in this case we have to remark that the uniform
convergence of the obstacles implies their level set convergence (via
remark~\rf{cap}). But the result we have obtained in this section does
not require that the obstacles be equicontrolled, and the convergence
of the solutions is in a different norm.}

\bigskip

{\bf Acknowldgement.} The author would like to thank prof.~Gianni Dal Maso
for helpful discussions and valuable suggestions.

\intro{References}
\def\interrefspace{\smallskip}  
{\ninepoint\frenchspacing

\item{[\bib{ATT}]}
ATTOUCH  H.:
{\it Variational convergence for 
functions and operators}. Pitman, Boston, 1984

\item{[\bib{ATT-PIC-1}]}
ATTOUCH  H., PICARD  C.:
Probl\`emes variationnels et th\'eorie du potentiel non lin\'eaire.
{\it Ann. Fac. Sci. Toulouse } {\bf 1} (1979), 89-136
\interrefspace 

\item{[\bib{ATT-PIC-2}]}
ATTOUCH  H., PICARD  C.:
In\'equations variationnelles avec obstacles et espaces fonctionnels en
th\'eorie du potentiel.
{\it Applicable Analysis } {\bf 12} (1981), 287-306
\interrefspace 


\item{[\bib{BOC-CIR-1}]}
BOCCARDO  L., CIRMI  G.R.:
Nonsmooth unilateral problems. {\it Nonsmooth optimization: methods and
applications (Erice, 1991)}, F. Giannessi ed., 1-10. Gordon and Breach,
Amsterdam, 1992
$L^1$-data.
\interrefspace

\item{[\bib{BOC-CIR-2}]}
BOCCARDO  L., CIRMI  G.R.:
Existence and uniqueness of solution of unilateral problems with
$L^1$-data.
\interrefspace

\item{[\bib{BOC-GAL}]}
BOCCARDO  L., GALLOU\"ET  T.:
Probl\`emes unilat\'eraux avec donn\'ees dans $\hbox{L}^1$. {\it
C. R. Acad. Sci. Paris S\'er. I Math.} {\bf 311} (1990), 617--655.
\interrefspace

\item{[\bib{BOC-GAL-ORS}]}BOCCARDO  L., GALLOU\"ET  T., ORSINA  L.:
Existence and uniqueness of entropy solutions for nonlinear elliptic
equations with measure data. {\it Ann. Inst. Henri Poincar\'e}  {\bf 13}
(1996), 539-551.
\interrefspace

\item{[\bib{BOC-MUR-1}]}
BOCCARDO  L., MURAT  F.:
Nouveaux r\'esultats de convergence dans des probl\`emes unilateraux.
In ``Nonlinear partial differential equations and their applications''.
Coll\`ege de France Seminar, Volume II. H. Br\'ezis and J.L. Lions
editors. Pitman, London, 1982
\interrefspace

\item{[\bib{DAL-3}]} 
DAL MASO G.: 
On the integral representation of certain local functionals. 
{\it Ricerche Mat.}, {\bf 22} (1983), 85-113.
\interrefspace

\item{[\bib{DAL-4}]}
DAL MASO  G.:
Some necessary and sufficient conditions for the convergence of sequences
of unilateral convex sets.
{\it J. Funct. Anal.} {\bf 62,2} (1985), 119-159.
\interrefspace

\item{[\bib{DAL-DAL}]}
DALL'AGLIO  P., DAL MASO G.:
Some properties of solutions of obstacle problems with measure data.
Preprint S.I.S.S.A., Trieste, 1998
\interrefspace

\item{[\bib{DAL-LEO}]}
DALL'AGLIO  P., LEONE  C.:
Obstacle problems with measure data.
Preprint S.I.S.S.A., Trieste, 1997.
\interrefspace

\item{[\bib{FUG}]}
FUGLEDE  B.:
The quasi topology associated with a countably subadditive set function.
{\it Ann. Inst. Fourier} {\bf 21,1} (1971), 123-169
\interrefspace

\item{[\bib{FUK}]}
FUKUSHIMA M., SATO K., TANIGUCHI S.:
On the closable part of pre-Dirichlet forms and the fine supports of 
underlying measures.
{\it Osaka J. Math.} {\bf 28} (1991), 517-535.
\interrefspace


\item{[\bib{HEI-KIL}]}HEINONEN  J., KILPEL\"AINEN  T., MARTIO  O.: {\it
Nonlinear potential theory of degenerate elliptic equations.} Clarendon
Press, Oxford, 1993.
\interrefspace

\item{[\bib{KIN-STA}]}KINDERLEHRER  D., STAMPACCHIA  G.: {\it
An introduction to variational inequalities and their applications.}
Academic, New York, 1980.


\item{[\bib{LEO}]}
LEONE  C.:
Existence and uniqueness of solutions for nonlinear obstacle problems
with
measure data.
Preprint S.I.S.S.A., Trieste, 1998
\interrefspace



\item{[\bib{MOS}]}
MOSCO  U.:
Convergence of convex sets and of solutions of variational inequalities. 
{\it Adv. in Math.\/} {\bf 3} (1969), 510-585
\interrefspace

\item{[\bib{OPP-ROS-1}]}
OPPEZZI P., ROSSI A.M.:
Existence and solutions for unilateral problems with multivalued
operators. {\it J. Convex. Anal.} {\bf 2} (1995) 241--261.
\interrefspace

\item{[\bib{OPP-ROS-2}]}
OPPEZZI P., ROSSI A.M.:
Esistenza di soluzioni per problemi unilaterali con dato misura o
$\hbox{L}^1$. {\it Ricerche Mat.}, to appear.
\interrefspace

\item{[\bib{ORS}]}
ORSINA  L.:
Solvability of linear and semilinear eigenvalue problems with $\Luno$
data.
{\it Rend. Sem. Mat. Padova } {\bf 90} (1993), 207-238
\interrefspace


\item{[\bib{STA}]}
STAMPACCHIA  G.:
Le probl\`eme de Dirichlet pour les \'equations elliptiques du second
ordre \`a coefficients discontinus.
{\it Ann. Inst. Fourier Grenoble\/} {\bf 15} (1965), 189-258.
\interrefspace

\item{[\bib{TRO}]}TROIANIELLO  G.M.: {\it
Elliptic differential equations and obstacle problems.}
Plenum Press, New York, 1987.

\end